\documentclass[11pt]{article}

\usepackage{amssymb,latexsym,amsmath}

\usepackage{cite}

\usepackage{eucal}

\title{Traces for star products on the dual of a Lie algebra}

\author{
    \textbf{
    Pierre Bieliavsky$^a$\thanks{pbiel@ulb.ac.be},
    \addtocounter{footnote}{1}
    Martin Bordemann$^b$\thanks{M.Bordemann@univ-mulhouse.fr},
    Simone Gutt$^a$\thanks{sgutt@ulb.ac.be},
    Stefan Waldmann$^c$\thanks{Stefan.Waldmann@physik.uni-freiburg.de}
    }
    \\[1cm]
    \begin{minipage} {7cm}
        \begin{center}
            $^a$D{\'e}partement de Math{\'e}matique \\
            Universit{\'e} Libre de Bruxelles \\
            Campus Plaine, C. P. 218 \\
            Boulevard du Triomphe \\
            B-1050 Bruxelles \\
            Belgique
        \end{center}
    \end{minipage}
    \begin{minipage} {7cm}
        \begin{center}
            $^b$Laboratoire de Math{\'e}matiques\\
            Universit{\'e} de Haute-Alsace Mulhouse \\
            4, Rue des Fr{\`e}res Lumi{\`e}re \\
            F.68093 Mulhouse \\
            France
        \end{center}
    \end{minipage}
    \\[2cm]
    \begin{minipage} {7cm}
        \begin{center}
            $^c$Fakult{\"a}t f{\"u}r Physik \\
            Albert-Ludwigs-Universit{\"a}t Freiburg \\
            Hermann Herder Stra{\ss}e 3 \\
            D 79104 Freiburg \\
            Germany
        \end{center}
    \end{minipage}
    }

\date{February 2002}

\newcommand{\im} {{\mathrm{i}}}
\newcommand{\eu} {{\mathrm{e}}}

\newcommand{\supp}       {\mathop{{\mathrm{supp}}}}
\newcommand{\image}      {\mathop{{\mathrm{im}}}\nolimits}
\newcommand{\ad}         {\mathop{{\mathrm{ad}}}\nolimits}

\newcommand{\End}        {\mathop{{\mathsf{End}}}\nolimits}

\newcommand{\tr}         {\mathop{{\mathsf{tr}}}\nolimits}
\newcommand{\Pol}        {\mathop{{\mathsf{Pol}}}\nolimits}

\newcommand{\ring}[1]    {{\mathsf {{#1}}}}

\newcommand{\cc}[1]      {\overline{{#1}}}

\newcommand{\SP}[1]      {{\left\langle{{#1}}\right\rangle}}

\newcommand{\Lie}        {\mathcal{L}}

\newcommand{\id}         {\mathsf{id}}

\newcommand{\prol}       {\mathrm{prol}}

\newcommand{\Left}        {\mathsf{L}}
\newcommand{\Right}      {\mathsf{R}}

\newcommand{\ind} {\mathrm{ind}}

\newcommand{\starWeyl} {\mathbin{*_{\scriptscriptstyle\mathrm{Weyl}}}}
\newcommand{\starG}    {\mathbin{*_{\scriptscriptstyle\mathrm{G}}}}
\newcommand{\starBCH}  {\mathbin{*_{\scriptscriptstyle\mathrm{BCH}}}}
\newcommand{\starK}    {\mathbin{*_{\scriptscriptstyle\mathrm{K}}}}
\newcommand{\starS}    {\mathbin{*_{\scriptscriptstyle\mathrm{S}}}}
\newcommand{\starW}    {\mathbin{*_{\scriptscriptstyle\mathrm{W}}}}
\newcommand{\starF}    {\mathbin{*_{\scriptscriptstyle\mathrm{F}}}}
\newcommand{\starO}    {\mathbin{*_{\scriptscriptstyle\mathcal{O}}}}
\newcommand{\starBM}   {\mathbin{*_{\scriptscriptstyle\mathrm{BM}}}}

\newtheorem{lemma}{Lemma}[section]
\newtheorem{proposition}[lemma]{Proposition}
\newtheorem{theorem}[lemma]{Theorem}
\newtheorem{corollary}[lemma]{Corollary}

\newtheorem{example}[lemma]{Example}
\newtheorem{remark}[lemma]{Remark}
\newtheorem{question}[lemma]{Question}

\newenvironment{proof}[1][{}]{\small{\textsc{Proof{#1}:~}}}{{\hspace*{\fill}$\square$\\}}

\numberwithin{equation}{section}

\begin{document}

\maketitle

\begin{abstract}
In this paper, we describe all traces for the BCH star-product on the dual
of a Lie algebra. First we show by an elementary argument 
that the BCH as well as the Kontsevich 
star-product are strongly closed if and only if the Lie algebra is 
unimodular. In a next step we show that the traces of the BCH 
star-product are given by the $\ad$-invariant functionals. 
Particular examples are the integration over coadjoint orbits. 
We show that for a compact Lie group and a regular orbit one can even 
achieve that this integration becomes a positive trace functional. 
In this case we explicitly describe the corresponding GNS representation.
Finally we discuss how invariant deformations on a group can be used to induce 
deformations of spaces where the group acts on.
\end{abstract}

%
%

\section{Introduction}

Trace functionals play an important role in deformation
quantization \cite{bayen.et.al:1978} (for recent reviews on
deformation quantization we refer to
\cite{gutt:2000,weinstein:1994,sternheimer:1998,dito.sternheimer:2002},
existence and classification results can be found in
\cite{kontsevich:1997:pre,fedosov:1996,omori.maeda.yoshioka:1991,nest.tsygan:1995a,bertelson.cahen.gutt:1997,weinstein.xu:1998}).

Physically, traces correspond to states of thermodynamical equilibrium
characterized by the KMS condition at infinite temperature
\cite{basart.et.al:1984,bordemann.roemer.waldmann:1998}. Note however,
that for reasonable physical interpretation one has to impose an
additional positivity condition on the traces
\cite{bordemann.waldmann:1998,waldmann:2000}.

On the mathematical side traces are one half of the index theorem,
namely the part of cyclic cohomology. The other half comes from the
$K$-theory part. Having a trace functional 
$\tr: \mathcal{A} \to \ring{C}$ of an associative algebra
$\mathcal{A}$ over some commutative ring $\ring{C}$ and having a
projection $P = P^2 \in M_n(\mathcal{A})$ representing an element 
$[P] \in K_0 (\mathcal{A})$ the value $\tr(P) \in \ring{C}$ does not
depend on $P$ but only on its class $[P]$. This is just the usual
natural pairing of cyclic cohomology with $K$-theory, see
e.g.~\cite[Chap.~III.3]{connes:1994}, and the value 
$\ind([P]) = \tr(P)$ is called the index of $[P]$ with respect to the
chosen trace.

In the case of deformation quantization the situation is as
follows. The starting point is a star-product $\star$ for a Poisson
manifold $(M, \pi)$ whence the algebra of interest is
$\mathcal{A} = (C^\infty(M)[[\nu]], \star)$ viewed as an algebra over
$\mathbb{C}[[\nu]]$. Then a trace is a $\mathbb{C}[[\nu]]$-linear
functional $\tr: C^\infty(M)[[\nu]] \to \mathbb{C}[[\nu]]$ such that
\begin{equation}
    \label{eq:atraceisatrace}
    \tr(f \star g) = \tr(g \star f),
\end{equation}
whenever one function has compact support. For the $K$-theory part of
the index theorem one knows that $K$-theory is stable under
deformation, see e.g.~\cite{rosenberg:1996:pre}: any projection $P_0$
of the undeformed algebra $M_n(C^\infty(M))$ can be deformed into a
projection
\begin{equation}
    \label{eq:fedosovsprojection}
    P = \frac{1}{2} + \left(P_0 - \frac{1}{2}\right) \star
    \frac{1}{\sqrt[\star]{1 + 4 (P_0 \star P_0 - P_0)}}
\end{equation}
with respect to $\star$, see
\cite[Eq.~(6.1.4)]{fedosov:1996}. Moreover, this deformation is unique
up to equivalence of projections and any projection of the deformed
algebra arises this way. It follows that $\ind([P])$ only depends on
$[P_0] \in K_0 (C^\infty(M))$, which is the isomorphism class of the
vector bundle defined by $P_0$, see also
\cite{bursztyn.waldmann:2000b} for a more detailed discussion.

Now let $\tilde{\star}$ be an equivalent star product with equivalence
transformation $T(f \star g) = Tf \tilde{\star} Tg$. Then clearly
$\tilde{\tr} = \tr \circ T^{-1}$ defines a trace functional with
respect to $\tilde{\star}$. From \eqref{eq:fedosovsprojection} we see
that $\ind([P]) = \tilde{\ind}([\tilde P])$ where $\tilde{\ind}$ is
the index with respect to the trace $\tilde{\tr}$ and
$\tilde{\star}$. Thus the index transforms well under equivalences of
star products provided one uses the `correct` corresponding trace. It
happens that in the symplectic case there is only one trace up to
normalization \cite{nest.tsygan:1995a}. So suppose that $M$ is compact
and that for each star product $\star$ we have chosen a trace
$\tr_\star$ normalized such that $$\tr_\star (1) = c$$ where $c$ does
not depend on $\star$. Then $T1 = 1$ implies 
$\tr_{\tilde{\star}} = \tr_\star \circ T^{-1}$ and thus the index does
not depend on the choice of $\star$ but only on the equivalence class
$[\star]$. This simple reasoning already explains the structure of
Fedosov's index formula \cite[Thm.~6.1.6]{fedosov:1996},
\cite{nest.tsygan:1995a}. Nevertheless we would like to mention that
the computation of $\ind([P])$ in geometrical terms is a quite
non-trivial task.

For a formulation of the index theorem in the general Poisson case we
refer to \cite{tamarkin.tsygan:2001}. Here the situation is far more
non-trivial as in general there is no longer a unique trace. In
\cite{felder.shoikhet:2000} it is shown that integration over $M$ with
respect to some smooth density $\Omega$ is a trace for Kontsevich's
star product provided the Poisson tensor is $\Omega$-divergence
free. However, there are much more traces, typically involving
integrations over the symplectic leaves.

An elementary proof that in the symplectic case 
one has a unique trace is presented in
\cite{gutt.rawnsley:2001:pre}. This approach
uses the canonical way of normalization of the trace,
introduced by Karabegov \cite{karabegov:1998b} using local $\nu$-Euler
derivations, see \cite{gutt.rawnsley:1999} and the elementary
proof of the uniqueness up to scaling of a trace as given in
\cite{bordemann.roemer.waldmann:1998}: Here one uses the fact that 
in the whole algebraic dual of $C^\infty(M)$ there is only one 
\emph{Poisson trace}
\begin{equation}
    \label{eq:PoissonTrace}
    \tau_0 (\{f, g\}) = 0,
\end{equation}
namely the integration with respect to the Liouville measure.

In this article we shall now consider the most simple case of 
a Poisson manifold: the dual of a Lie algebra. Here we shall 
determine all the traces for the BCH star product on 
$\mathfrak{g}^*$ by very elementary arguments.

The paper is organized as follows. In
Sect.~\ref{sec:notation} we recall the construction of various star
products on the dual of a Lie algebra $\mathfrak{g}^*$ as well as
their relation to star products on $T^*G$ where $G$ is a Lie group
with Lie algebra $\mathfrak{g}$. Then we prove the strong closedness
of homogeneous star products on $\mathfrak{g}^*$ by elementary
computations in Sect.~\ref{sec:closed} and in
Sect.~\ref{sec:invariant} we show that any $\ad$-invariant functional
is a trace for the BCH star product. In Sect.~\ref{sec:positive} we
prove the positivity of a trace $\boldsymbol{\tau_{\mathcal{O}}}$
associated to a regular orbit $\mathcal{O} \subseteq \mathfrak{g}^*$
for compact $G$ by a BRST construction of a star product on
$\mathcal{O}$. Sect.~\ref{sec:gns} contains a characterization of the
GNS representation obtained from the positive trace
$\boldsymbol{\tau_{\mathcal{O}}}$. Finally, Sect.~\ref{sec:action} is
devoted to a construction of trace functionals by a group action using
a `universal deformation' on the group, inspired by techniques
developed in \cite{giaquinto.zhang:1998,bieliavsky.massar:2000:pre}.

%
%

\section{Star products on $\mathfrak{g}^*$ and $T^*G$}
\label{sec:notation}

In this section we shall recall the construction of several star
products on the dual $\mathfrak{g}^*$ of a Lie algebra $\mathfrak{g}$
and on $T^*G$ where $G$ is a Lie group with Lie algebra
$\mathfrak{g}$. First we shall establish some notation.

By $e_1, \ldots, e_n$ we denote a basis of $\mathfrak{g}$ with dual
basis $e^1, \ldots, e^n \in \mathfrak{g}^*$. Such a basis gives raise
to linear coordinates $x = x^i e_i$ on $\mathfrak{g}$ and
$\xi = \xi_i e^i$ on $\mathfrak{g}^*$. Here and in the following we
shall use Einstein's summation convention. With a capital letter $X$
we shall denote the left-invariant vector field
$X \in \Gamma^\infty(TG)$ corresponding to $x \in \mathfrak{g}$,
i.e. $X_e = x$. A vector $x \in \mathfrak{g}$ determines a linear
function $\hat x \in \Pol^1(\mathfrak{g}^*)$ by
$\hat x (\xi) = \xi(x)$. Analogously, $X \in \Gamma^\infty(TG)$
determines a function $\hat X \in \Pol^1(T^*G)$, linear in the
fibers, by setting $\hat X(\alpha_g) = \alpha_g (X_g)$, where
$\alpha_g \in T^*_g G$ and $g \in G$. We shall use the same symbol
$\;\hat{}\;$ for the corresponding graded algebra isomorphism between
the symmetric algebra $\bigvee^\bullet\mathfrak{g}$ of $\mathfrak{g}$
and all polynomials $\Pol^\bullet(\mathfrak{g}^*)$ on
$\mathfrak{g}^*$. Similar we have a graded algebra isomorphism between
$\Gamma(\bigvee^\bullet TG)$ and $\Pol^\bullet(T^*G)$. By use of
left-invariant vector fields and one-forms, $TG$ and $T^*G$ trivialize
canonically. This yields $TG \cong G \times \mathfrak{g}$ and
$T^*G \cong G \times \mathfrak{g}^*$. The corresponding projections are
denoted by
\begin{equation}
    \label{eq:PiRhoProj}
    G \stackrel{\pi}{\longleftarrow} G \times \mathfrak{g}^*
    \stackrel{\varrho}{\longrightarrow} \mathfrak{g}^*,
\end{equation}
whence in particular $\hat X = \varrho^* \hat x$ for a left-invariant
vector field $X$. More generally,
$\Pol^\bullet (T^*G)^G = \varrho^* \Pol^\bullet(\mathfrak{g}^*)$. For
the symplectic Poisson bracket on $T^*G$ we use the sign convention
such that the map $\hat{}\, : \Gamma^\infty (TG) \to \Pol^1(T^*G)$
becomes an isomorphism of Lie algebras (and not an anti-isomorphism as
in \cite{bordemann.neumaier.waldmann:1998}). Then the canonical linear
Poisson bracket on $\mathfrak{g}^*$ can be obtained by the observation
that left-invariant functions on $T^*G$ (with respect to the lifted
action) are a Poisson sub-algebra which is in linear bijection with
$C^\infty(\mathfrak{g}^*)$ via $\varrho^*$. Thus it is meaningful to
require $\varrho^*$ to be a morphism of Poisson algebras. In the global
coordinates $\xi_1,\ldots, \xi_n$ the resulting Poisson bracket on
$\mathfrak{g}^*$ reads as
\begin{equation}
    \label{eq:PoissonBracketg}
    \{f, g\} = \xi_k c^k_{ij}
    \frac{\partial f}{\partial \xi_i}
    \frac{\partial g}{\partial \xi_j},
\end{equation}
where $c^k_{ij} = e^k ([e_i, e_j])$ are the structure constants of
$\mathfrak{g}$ and $f, g \in C^\infty (\mathfrak{g}^*)$.

The first star-product on $\mathfrak{g}^*$ is essentially given by the
Baker-Campbell-Hausdorff series of $\mathfrak{g}$. One uses the total
symmetrization map
$\sigma_\nu: \Pol^\bullet (\mathfrak{g}^*)[\nu]
\to \mathcal U (\mathfrak{g})[\nu]$ into the universal enveloping
algebra of $\mathfrak{g}$, defined by
\begin{equation}
    \label{eq:SigmaNuDef}
    \sigma_\nu (\widehat{x_1 \vee \cdots \vee x_k})
    = \frac{\nu^k}{k!} \sum_{\tau \in S_k}
    x_{\tau(1)} \bullet \cdots \bullet x_{\tau(k)},
\end{equation}
where we have built in the formal parameter $\nu$ already at this
stage. Then
\begin{equation}
    \label{eq:StarBCHDef}
    \sigma_\nu (f \starBCH g) =
    \sigma_\nu(f) \bullet \sigma_\nu (g)
\end{equation}
yields indeed a deformed product $\starBCH$ for
$f, g \in \Pol^\bullet (\mathfrak{g}^*)[\nu]$, which turns out to
extend to a differential star-product for
$C^\infty (\mathfrak{g}^*)[[\nu]]$, see \cite{gutt:1983} for a detailed
discussion. Here we shall just mention a few properties of
$\starBCH$. First, $\starBCH$ is
\emph{strongly $\mathfrak{g}$-invariant}, i.e. for
$f \in C^\infty (\mathfrak{g}^*)[[\lambda]]$ and $x \in \mathfrak{g}$ we
have
\begin{equation}
    \label{eq:BCHStrongInv}
    \hat x \starBCH f - f \starBCH \hat x = \nu \{\hat x, f\}.
\end{equation}
Moreover, $\starBCH$ is \emph{homogeneous}: this means that the
operator
$\mathcal{H} = \nu \frac{\partial}{\partial \nu} + \Lie_E$,
where $E = \xi_i \frac{\partial}{\partial \xi_i}$ is the
\emph{Euler vector field}, is a derivation of $\starBCH$, i.e.
\begin{equation}
    \label{eq:BCHHomogeneous}
    \mathcal{H}(f \starBCH g) =
    \mathcal{H} f \starBCH g + f \starBCH \mathcal{H} g
\end{equation}
for all $f, g \in C^\infty(\mathfrak{g}^*)[[\nu]]$. It follows
immediately that $\Pol^\bullet (\mathfrak{g}^*)[\nu]$ is a
`convergent' sub-algebra generated by the constant and linear
polynomials. The relation to the BCH series can be seen as follows:
Consider the exponential functions $e_x (\xi) := \eu^{\xi(x)}$. Then
for all $x, y \in \mathfrak{g}$ one has
\begin{equation}
    \label{eq:ExpBCH}
    e_x \starBCH e_y = e_{\frac{1}{\nu} H(\nu x, \nu y)},
\end{equation}
where $H(\cdot, \cdot)$ is the BCH series of
$\mathfrak{g}$. Since bidifferential operators on $\mathfrak{g}^*$ are
already determined by their values on the exponential functions
$e_x$, $x \in \mathfrak{g}$, the star-product $\starBCH$ is already
determined by (\ref{eq:ExpBCH}). For a more detailed analysis and proofs
of the above statements we refer to
\cite{gutt:1983,bordemann.neumaier.waldmann:1998}.

The other star-product we shall mention is the Kontsevich star-product
$\starK$ for $\mathfrak{g}^*$. His general construction of a star
product for arbitrary Poisson structures on $\mathbb R^n$ simplifies
drastically in the case of a linear Poisson structure
(\ref{eq:PoissonBracketg}). We shall not enter the general construction
but refer to
\cite{kontsevich:1997:pre,kontsevich:1999,arnal.benamar.masmoudi:1999,arnal.benamar:2000,dito:1999}
for more details and just mention a few properties of $\starK$. First,
$\starK$ is \emph{$\mathfrak{g}$-covariant}, i.e. one has
\begin{equation}
    \label{eq:KCovariant}
    \hat x \starK \hat y - \hat y \starK \hat x
    = \nu \{\hat x, \hat y\}
    = \nu \widehat{[x,y]}
\end{equation}
for all $x, y \in \mathfrak{g}$. This is a weaker compatibility with
the (classical) $\mathfrak{g}$-action than (\ref{eq:BCHStrongInv}).
Moreover, $\starK$ is homogeneous, too, but in general $\starK$ and
$\starBCH$ do \emph{not} coincide but are only equivalent, see
\cite{dito:1999}.

Let us now recall how the star-product $\starBCH$ on $\mathfrak{g}^*$
is related to star-products on $T^*G$. The main idea is to make the
Poisson morphism $\varrho^*$ into an algebra morphism of star-product
algebras. This requirement does not determine the star-product on
$T^*G$ completely and the remaining freedom (essentially the choice of
an `ordering prescription' between functions depending only on $G$ and
on $\mathfrak{g}^*$, respectively) can be used to impose further
properties. In \cite{gutt:1983} a star-product $\starG$ of Weyl-type was
constructed by inserting additional derivatives in $G$-direction into
the bidifferential operators of $\starBCH$. In
\cite{bordemann.neumaier.waldmann:1998} a star
product $\starS$ of standard-ordered type was obtained by a
(standard-ordered) Fedosov construction using the lift of the
half-commutator connection on $G$ to a symplectic connection on
$T^*G$. The star-product $\starS$ can also be understood as the
resulting composition law of symbols from the standard-ordered symbol
and differential operator calculus  induced by the half-commutator
connection. A further `Weyl-symmetrization' yields a star-product
$\starW$ of Weyl-type which does not coincide in general with the
original Fedosov star-product $\starF$ built out of the
half-commutator connection directly. However, it was shown in
\cite[Sect~8]{bordemann.neumaier.waldmann:1998} that $\starW$
coincides with $\starG$. Moreover, the pull-back $\varrho^*$ is indeed
an algebra morphism for both star-products $\starG$ and $\starS$,
i.e. one has
\begin{equation}
    \label{eq:PullBack}
    \varrho^*f *_{\scriptscriptstyle\mathrm{G/S}} \varrho^* g =
    \varrho^* (f \starBCH g)
\end{equation}
for all $f, g \in C^\infty (\mathfrak{g}^*)[[\nu]]]$. All the star
products $\starG$, $\starS$, and $\starF$ are homogeneous in the sense
of star-products on cotangent bundles whence it follows that they are
all strongly closed: integration over $T^*G$ with respect to the
Liouville form defines a trace on the functions with compact support,
see \cite[Sect.~8]{bordemann.neumaier.waldmann:1999}.

%
%

\section{Strong closedness of $\starBCH$ and $\starK$}
\label{sec:closed}

We shall now discuss an elementary proof of the fact that $\starBCH$
as well as $\starK$ are strongly closed with respect to the constant
volume form $d^n\xi$ on $\mathfrak{g}^*$ if and only if the Lie algebra
$\mathfrak{g}$ is \emph{unimodular}, i.e.
$\tr\ad(x) = 0$ for all $x \in \mathfrak{g}$, or, equivalently,
$c^i_{ij} = 0$. The unimodularity of $\mathfrak{g}$ is easily seen to
be necessary since it is exactly the condition that the integration is
a \emph{Poisson trace}, see also \cite[Sect.~4]{weinstein:1997} for
the Poisson case and \cite{felder.shoikhet:2000} for a different and
more general proof for Kontsevich's star product on $\mathbb{R}^n$.

Before we discuss the general case let us consider the case where $G$
is compact. In this case $\mathfrak{g}$ is known to be in particular
unimodular.
\begin{proposition}
    \label{proposition:CompactCase}
    Let $G$ be compact. Then $\starBCH$ is strongly closed.
\end{proposition}
\begin{proof}
    Let $f, g \in C^\infty_0 (\mathfrak{g}^*)$. Since $G$ is compact,
    $\varrho^*f, \varrho^*g \in C^\infty_0 (T^*G)$ and thus the strong
    closedness of $\starG$ and (\ref{eq:PullBack}) implies
    \[
    0 =
    \int_{T^*G} (\varrho^*f \starG \varrho^*g
    - \varrho^*g \starG \varrho^*f) \; \Omega
    = \mathrm{vol} (G) \int_{\mathfrak{g}^*}
    (f \starBCH g - g \starBCH f) \; d^n\xi,
    \]
    where $\Omega$ is the (suitably normalized) Liouville measure on
    $T^*G$.
\end{proof}

Clearly the above proof relies on the compactness of $G$, otherwise
the integration would not be defined. As an amusing observation we
remark that one can also use the above proposition to obtain the
well-known fact that compact Lie groups have unimodular Lie algebras.

For the general unimodular case we use a different argument which is
essentially the same as for homogeneous star-products on a cotangent
bundle \cite[Sect.~8]{bordemann.neumaier.waldmann:1999}. A
differential operator $D$ on $\mathfrak{g}^*$ is called homogeneous of
degree $r \in \mathbb Z$ if $[\Lie_E, D] = rD$, where $\Lie_E$ is the
Lie derivative with respect to the Euler vector field.
\begin{lemma}
    \label{lemma:Diffop}
    Let $D$ be a homogeneous differential operator of degree $-r$ with
    $r \ge 1$. Then for all $f \in C^\infty_0 (\mathfrak{g}^*)$ one has
    \begin{equation}
        \label{eq:IntHomDiffNull}
        \int_{\mathfrak{g}^*} Df \; d^n\xi = 0.
    \end{equation}
\end{lemma}
From here we can follow \cite{bordemann.neumaier.waldmann:1999} almost
literally: If $f \in \Pol^k (\mathfrak{g}^*)$ and
$g \in C^\infty_0 (\mathfrak{g}^*)$ then for every homogeneous star
product $*$ on $\mathfrak{g}^*$ one has
\begin{equation}
    \label{eq:IntPol}
    \int_{\mathfrak{g}^*} f * g \; d^n \xi =
    \sum_{r=0}^k \nu^r
    \int_{\mathfrak{g}^*} C_r(f, g) \; d^n \xi,
\end{equation}
where $C_r$ is the $r$-th bidifferential operator of $*$. This
follows from Lemma~\ref{lemma:Diffop} since $C_r (f, \cdot)$ is
homogeneous of degree $k-r$. The analogous statement holds for the
integral over $g*f$. From this we conclude the following lemma:
\begin{lemma}
    \label{lemma:PolTrace}
    Let $*$ be a homogeneous star-product for $\mathfrak{g}^*$,
    $f \in \Pol^\bullet (\mathfrak{g}^*)$, and
    $g \in C^\infty_0 (\mathfrak{g}^*)$. Then
    \begin{equation}
        \label{eq:TracePol}
        \int_{\mathfrak{g}^*} (f *g - g * f) \; d^n \xi = 0
    \end{equation}
    if and only if $\mathfrak{g}$ is unimodular.
\end{lemma}
\begin{proof}
    The proof is done by induction on the polynomial degree $k$ of
    $f$. For $k = 0$ the statement (\ref{eq:TracePol}) is true by
    (\ref{eq:IntPol}). For $k = 1$ we obtain (\ref{eq:TracePol}) by
    (\ref{eq:IntPol}) if and only if the integral vanishes on Poisson
    brackets, i.e. if and only if $\mathfrak{g}$ is unimodular. For
    $k \ge 2$ we can write $f$ as a $*$-polynomial in at most linear
    polynomials since these polynomials generate
    $\Pol^\bullet (\mathfrak{g}^*)[\nu]$ by the homogeneity of
    $*$. Then we can use the cases $k=0,1$ to prove (\ref{eq:TracePol}).
\end{proof}

Having the trace property for polynomials and compactly supported
functions, we only have to use a density argument, i.e. the
Stone-Weierstra{\ss} theorem, to conclude the trace property in
general:
\begin{theorem}
    \label{theorem:StrongClosedness}
    Let $*$ be a homogeneous star-product for $\mathfrak{g}^*$. Then
    the integration over $\mathfrak{g}^*$ with respect to the constant
    volume $d^n\xi$ is a trace if and only if $\mathfrak{g}$ is
    unimodular.
\end{theorem}
Since $\starBCH$ as well as $\starK$ are homogeneous this theorem
proves in an elementary way that they are strongly closed in the sense
of \cite{connes.flato.sternheimer:1992}.

%
%

\section{Trace properties of $\mathfrak{g}$-invariant functionals}
\label{sec:invariant}

Quite contrary to the symplectic case it turns out that in the Poisson
case traces are no longer unique in general.

Before we give an elementary proof in the case of $\mathfrak{g}^*$ we
shall make a few comments on the general situation. As we have seen
already before, the trace functionals are typically not defined on the
whole algebra but on a certain subspace, as e.g. the functions with
compact support. On the other hand, the property of being a trace only
becomes interesting if this subspace is not only a sub-algebra but even
an ideal. This motivates the following terminology:
For an associative algebra $\mathcal A$ we call a functional $\tau$
defined on $\mathcal{J} \subseteq \mathcal A$ a
\emph{trace on $\mathcal{J}$}
if $\mathcal{J}$ is a two-sided ideal and for all $A \in \mathcal{A}$
and $B \in \mathcal{J}$ one has $\tau([A, B]) = 0$. Similarly we define
a \emph{Poisson trace on a Poisson ideal} of a Poisson algebra.

With this notation the traces which are given by integrations are
traces on the ideals $C^\infty_0 (\mathfrak{g}^*)$
and $C^\infty_0(\mathfrak{g}^*)[[\nu]]$, respectively. However, there
will be some interesting traces with a slightly different domain. If
we want to integrate over a sub-manifold $\iota: N \hookrightarrow M$
then the following space becomes important. Here and in the following
we shall only consider the case where $\iota$ is an embedding. We
define
\begin{equation}
    \label{eq:CNMDef}
    C^\infty_N (M) := \{f \in C^\infty (M) \; | \;
                      \iota (N) \cap \supp f
                      \textrm{ is compact } \}.
\end{equation}
If $N$ is a closed embedded sub-manifold then
$C^\infty_0 (M) \subseteq C^\infty_N (M)$. Moreover, the locality of a
star-product ensures that $C^\infty_N (M)[[\nu]]$ is a two-sided ideal
of $C^\infty (M)[[\nu]]$.

Taking such a subspace as example we consider more generally
domains of the form $\mathcal{D}[[\nu]]$ where
$\mathcal{D} \subseteq C^\infty(M)$. In this case $\mathcal{D}$ is
necessarily a Poisson ideal which follows immediately from the ideal
properties of $\mathcal{D}[[\nu]]$. Moreover, if
$\tau: \mathcal{D}[[\nu]] \to \mathbb{R}[[\nu]]$ is a trace for a local
star-product $*$ on $M$ with domain $\mathcal{D}[[\nu]]$ then
$\tau = \sum_{r=0}^\infty \nu^r \tau_r$ with linear functionals
$\tau_r: \mathcal{D} \to \mathbb{R}$. For the following we shall assume
that all $\tau_r$ have some reasonable continuity property, e.g. with
respect to the locally convex topology of smooth functions. This
requirement seems to be reasonable as long as we are dealing with star
products having at least continuous cochains in every order of $\nu$.

Now let us come back to the case of $\mathfrak{g}^*$ with the star
product $\starBCH$. As a first observation we remark that the strong
$\mathfrak{g}$-invariance of $\starBCH$ implies that for a two-sided
ideal $\mathcal{D}[[\nu]]$ the space $\mathcal{D}$ is
$\mathfrak{g}$-invariant. Moreover, we have the following theorem:
\begin{theorem}
    \label{theorem:Trace}
    Let $\mathcal{D} \subseteq C^\infty (\mathfrak{g}^*)$ be a subspace
    such that $\mathcal{D}[[\nu]]$ is a two-sided ideal with respect to
    $\starBCH$ and let $\tau = \sum_{r=0}^\infty \nu^r \tau_r$ be a
    $\mathbb{R}[[\nu]]$-linear functional on $\mathcal{D}[[\nu]]$ with
    the following continuity property: For a given
    $f \in C^\infty (\mathfrak{g}^*)$ and $g \in \mathcal{D}$ and a
    sequence $p_n \in \Pol^\bullet (\mathfrak{g}^*)$ such that
    $p_n \to f$ in the locally convex topology of smooth functions we
    have $\tau_r ([p_n, g]_{\starBCH}) \to \tau_r ([f, g]_{\starBCH})$
    (in each order of $\nu$).

    Then $\tau$ is a $\starBCH$-trace on $\mathcal{D}[[\nu]]$ if and
    only if $\tau$ is a Poisson trace on $\mathcal{D}$ which is the
    case if and only if $\tau$ is $\mathfrak{g}$-invariant.
\end{theorem}
\begin{proof}
    The continuity ensures that $\mathfrak{g}$-invariance coincides with
    the property of being a Poisson trace. Now let $\tau_0$ be a
    Poisson trace and let $g \in \mathcal{D}$.
    Then for all $x \in \mathfrak{g}$ we have
    $\tau_0 ([\hat x, g]) = \nu \tau_0 (\{\hat x, g\}) = 0$ by the
    strong invariance of $\starBCH$. But since
    $\Pol^1(\mathfrak{g}^*)[\nu]$ together with the constants generates
    $\Pol^\bullet (\mathfrak{g}^*)[\nu]$ we have $\tau_0 ([p,g]) = 0$
    for every polynomial $p$. Together with the fact that the
    polynomials are dense in $C^\infty (\mathfrak{g}^*)$ and $\tau_0$
    has the above continuity it follows that $\tau_0$ is a
    $\starBCH$-trace.
    Now if $\tau$ is a $\starBCH$-trace then $\tau_0$ is a Poisson
    trace and hence a $\starBCH$-trace itself. Thus $\tau - \tau_0$ is
    still a $\starBCH$-trace and a simple induction proves the
    theorem.
\end{proof}

The somehow technical continuity property needed above turns out to be
rather mild. In the main example it is trivially fulfilled:
\begin{example}
    \label{example:OrbitIntegration}
    \hfill
    \begin{enumerate}
    \item Let $\iota: \mathcal{O} \hookrightarrow \mathfrak{g}^*$ be a
        not necessarily closed but embedded coadjoint orbit and
        consider
        $\mathcal{D} = C^\infty_{\mathcal{O}} (\mathfrak{g}^*)$. Then the
        integration with respect to the Liouville measure
        $\Omega_{\mathcal{O}}$ on $\mathcal{O}$,
        \begin{equation}
            \label{eq:TauODef}
            \tau_{\mathcal{O}} (f) := \int_{\mathcal{O}} \iota^* f
            \; \Omega_{\mathcal{O}},
        \end{equation}
        is a $\starBCH$-trace on
        $C^\infty_{\mathcal{O}} (\mathfrak{g}^*)[[\nu]]$.
    \item If in addition $\Delta$ is a $\mathfrak{g}$-invariant
        differential operator on $\mathfrak{g}^*$ then
        $\tau^\Delta_{\mathcal{O}}$, defined by
        \begin{equation}
            \label{eq:TauODeltaDef}
            \tau^\Delta_{\mathcal{O}} (f) :=
            \tau_{\mathcal{O}} (\Delta f) = \int_{\mathcal{O}}
\iota^*(\Delta f)
            \; \Omega_{\mathcal{O}},
        \end{equation}
        is still a trace on
        $C^\infty_{\mathcal{O}} (\mathfrak{g}^*)[[\nu]]$.
    \end{enumerate}
\end{example}

%
%

\section{Positivity of traces}
\label{sec:positive}

If one replaces the formal parameter $\nu$ by a new formal parameter
$\lambda$ such that $\nu = \im \lambda$ and if one treats $\lambda$ as
a real quantity, i.e. $\cc\lambda = \lambda$, then it is well-known
that the complex conjugation of functions in
$C^\infty(\mathfrak{g}^*)[[\lambda]]$ becomes a $^*$-involution for
$\starBCH$. One has
\begin{equation}
    \label{eq:CCInvBCH}
    \cc{f \starBCH g} = \cc{g} \starBCH \cc{f}
\end{equation}
for all $f, g \in C^\infty(\mathfrak{g})[[\lambda]]$. Such a star
product is also called a Hermitian star-product, see
e.g.~\cite{bursztyn.waldmann:2000a} for a detailed discussion. Thus one
enters the realm of $^*$-algebras over ordered rings, see
\cite{bursztyn.waldmann:2001a,bordemann.waldmann:1998}. In particular
one can ask whether the traces for $\starBCH$ are \emph{positive}
linear functionals, i.e. satisfy $\tau( \cc f \starBCH f) \ge 0$ in
the sense of formal power series, if the corresponding classical
functional $\tau_0$ comes from a positive Borel measure on
$\mathfrak{g}$. In general a classically positive linear functional is
no longer positive for a deformed product, see
e.g.~\cite[Sect.~2]{bordemann.waldmann:1998} for
a simple example and \cite{bursztyn.waldmann:2000a}. But sometimes one
can \emph{deform} the functional as well in order to make it positive
again: in the case of star-products on symplectic manifolds this is
always possible \cite[Prop.~5.1]{bursztyn.waldmann:2000a}. Such
deformations are called \emph{positive deformations}. In our case we
are faced with the question whether we can deform the traces
$\tau_{\mathcal{O}}$ such that on one hand they are still traces and
on the other hand they are positive.

One strategy could be the following: First prove that the trace can be
deformed into a positive functional perhaps loosing the trace
property. Secondly average over the group in order to obtain a
$\mathfrak{g}$-invariant functional and hence a trace. This would
require to have a compact group. However, we shall follow another idea
giving some additional insight in the problem. Nevertheless we shall
first ask the following question as a general problem in deformation
quantization of Poisson manifolds:
\begin{question}
    \label{question:PosDef}
    Is every Hermitian star-product on a Poisson manifold a positive
    deformation?
\end{question}

We shall now consider the following more particular case. We assume
the group $G$ to be compact and
$\iota: \mathcal{O} \hookrightarrow \mathfrak{g}^*$ to be a regular
coadjoint orbit. Then we want to find a positive trace for $\starBCH$
with zeroth order given by $\tau_{\mathcal{O}}$ as in
(\ref{eq:TauODef}). The construction is based on the following theorem
which is of independent interest:
\begin{theorem}
    \label{theorem:iotaHomo}
    Let $G$ be compact and let
    $\iota: \mathcal{O} \hookrightarrow \mathfrak{g}^*$ be a regular
    coadjoint orbit. Then there exists a star-product $\starO$ on the
    symplectic manifold $\mathcal{O}$ and a series of
    $\mathfrak{g}$-invariant differential operators
    $S = \id + \sum_{r=1}^\infty \lambda^r S_r$ on $\mathfrak{g}^*$
    such that the deformed restriction map
    \begin{equation}
        \label{eq:DefRest}
        \boldsymbol{\iota^*} = \iota^* \circ S:
        C^\infty (\mathfrak{g}^*)[[\lambda]] \to
        C^\infty (\mathcal{O})[[\lambda]]
    \end{equation}
    becomes a real surjective homomorphism of star-products, i.e.
    \begin{equation}
        \label{eq:DefRestHom}
        \boldsymbol{\iota^*}f \starO \boldsymbol{\iota^*} g =
        \boldsymbol{\iota^*} (f \starBCH g)
        \quad
        \textrm{and}
        \quad
        \cc{(\boldsymbol{\iota^*}f)} = \boldsymbol{\iota^*} \cc f
    \end{equation}
    for all $f, g \in C^\infty(\mathfrak{g}^*)[[\lambda]]$. Hence
    $\starO$ becomes a Hermitian deformation.
\end{theorem}
One can view this theorem as a certain `deformed tangentiality
property' of the star product $\starBCH$: Though $\starBCH$ is not
tangential, i.e. restricts to all orbits, for a particular
orbit it can be arranged such that it restricts by deforming the
restriction map, see \cite{cahen.gutt.rawnsley:1996} for a more
detailed discussion.

From this theorem and
\cite[Lem.~2]{bordemann.waldmann:1998} we immediately obtain a
positive trace deforming $\tau_{\mathcal{O}}$:
\begin{corollary}
    \label{corollary:posTrace}
    Let $G$ be compact and
    $\iota: \mathcal{O} \hookrightarrow \mathfrak{g}^*$ a regular
    orbit with deformed restriction map $\boldsymbol{\iota^*}$ as in
    (\ref{eq:DefRest}). Then the functional
    \begin{equation}
        \label{eq:posTrace}
        \boldsymbol{\tau}_{\mathcal{O}} (f) := \int_{\mathcal{O}}
        \boldsymbol{\iota^*} f \; \Omega_{\mathcal{O}}
    \end{equation}
    is a positive trace with classical limit
    $\tau_{\mathcal{O}}$. In particular, $\starO$ is strongly closed.
\end{corollary}
Thus it remains to prove Theorem~\ref{theorem:iotaHomo}. We shall use
here arguments from phase space reduction of star-products via the
BRST formalism as discussed in detail in
\cite{bordemann.herbig.waldmann:2000}. In order to make this article
self-contained we shall recall the basic steps of
\cite{bordemann.herbig.waldmann:2000} adapted to the case of Poisson
manifolds.

\smallskip

\noindent
\begin{proof}[ of Theorem~\ref{theorem:iotaHomo}]
    Since $\mathcal{O}$ is assumed to be a regular orbit there are 
    real-valued \emph{Casimir polynomials}
    $J_1, \ldots, J_k \in \Pol^\bullet(\mathfrak g^*)$ such that
    $\mathcal{O}$ can be written as level surface
    $\mathcal{O} = J^{-1} (\{0\})$ for
    the map 
    $J = (J_1, \ldots, J_k): \mathfrak{g}^* \to \mathbb{R}^k$, where
    $0$ is a regular value. Since the components of $J$ commute with
    respect to the Poisson bracket this can be viewed as a moment map
    $J: \mathfrak{g}^* \to \mathfrak{t}^*$ where $\mathfrak{t}^*$ is
    the dual of the $k$-dimensional Abelian Lie algebra. Moreover, the
    $J$'s are in the Poisson center whence the corresponding torus
    action is trivial.

    Since the differential operators $S_r$ will only be needed near
    $\mathcal{O}$ it will be sufficient to construct them in a tubular
    neighbourhood around $\mathcal{O}$. In fact, a globalization
    beyond is also easily obtained, see
    \cite[Lem.~6]{bordemann.herbig.waldmann:2000}. As $0$ is a regular
    value of $J$ we can use $J$ for the transversal coordinates and
    find a $G$-invariant tubular neighbourhood $U$ of
    $\mathcal{O}$. On $U$ we can define the following maps: First we
    need a prolongation map 
    $\prol: C^\infty(\mathcal{O}) \hookrightarrow C^\infty (U)$ given
    by 
    \begin{equation}
        \label{eq:prol}
        (\prol \phi)(o, \mu) = \phi (o),
    \end{equation}
    where $o \in \mathcal{O}$ and $\mu \in \mathfrak{t}^*$ is the
    transversal coordinate in $U$. Next we consider
    $\bigwedge^\bullet(\mathfrak{t}) \otimes C^\infty(\mathfrak{g}^*)$
    and define the \emph{Koszul coboundary operator} $\partial$ by the
    (left-)insertion of $J$, i.e.
    $\partial(t \otimes f) = \sum_l i(e^l) t \otimes J_l f$,
    where $J = \sum_l e^l J_l$. Clearly $\partial$ is $G$-invariant with
    respect to the $G$ action $g^* (t \otimes f) = t \otimes g^* f$
    and $\partial^2 = 0$. We shall denote the homogeneous components of
    $\partial$ by
    $\partial_l:
    \bigwedge^l (\mathfrak{t}) \otimes C^\infty (\mathfrak{g}^*)
    \to
    \bigwedge^{l-1} (\mathfrak{t}) \otimes C^\infty (\mathfrak{g}^*)$
    for $l \ge 1$. In the case $l = 0$ we set $\partial_0 = \iota^*$
    and clearly $\iota^*\partial_1 = 0$. Finally, we define the chain
    homotopy $h$ on
    $\bigwedge^\bullet (\mathfrak{t}) \otimes C^\infty (U)$ by
    \begin{equation}
        \label{eq:homotopy}
        h(t \otimes f) (o, \mu)
        = \sum_{l=1}^k e_l \wedge t \otimes \int_0^1
        \frac{\partial f}{\partial \mu_l} (o, s\mu) s^k ds,
    \end{equation}
    an denote the corresponding homogeneous components by $h_l$. For
    convenience we set $h_{-1} = \prol$. Then $h$ is obviously
    $G$-invariant and it is indeed a chain homotopy for $\partial$,
    i.e. for all $l = 0, \ldots, k$ we have
    \begin{equation}
        \label{eq:ch}
        h_{l-1} \partial_l + \partial_{l+1} h_l =
        \id_{\bigwedge^l (\mathfrak{t}) \otimes C^\infty (U)}.
    \end{equation}
    Moreover, one has the obvious identities
    \begin{equation}
        \label{eq:iprol}
        \iota^* \prol = \id_{C^\infty (\mathcal{O})},
        \quad
        \textrm{and}
        \quad
        h_0 \prol = 0.
    \end{equation}

    In a next step we quantize the above chain complex and it's
    homotopy. The first easy observation is that the star-product
    $\starBCH$ is strongly $\mathfrak{t}$-invariant, i.e. the
    components of $J$ are in the center of $\starBCH$, too. Thus we
    can define a deformed Koszul operator $\boldsymbol{\partial}$ on
    the space 
    $\left(\bigwedge^\bullet(\mathfrak{t})
        \otimes C^\infty (\mathfrak{g}^*)\right)[[\lambda]]$ by
    \begin{equation}
        \label{eq:defkuszul}
        \boldsymbol{\partial} (t \otimes f) 
        = \sum_l i(e^l) t \otimes f \starBCH J_l.
    \end{equation}
    Then we still have $\boldsymbol{\partial}^2 = 0$ as well as
    $\cc{\boldsymbol{\partial} (t \otimes f)}
    = \boldsymbol{\partial} \cc{(t \otimes f)}$
    since the $J_l$ commute and are real. Moreover,
    $\boldsymbol{\partial}$ is still $G$-invariant. In a next step one
    constructs the deformations of $h$ and $\iota^*$ as follows. We define
    $\boldsymbol{h}_{-1} = \prol$ without deformation and set
    \begin{equation}
        \label{eq:defhomotopies}
        \boldsymbol{\partial}_0 := \boldsymbol{\iota^*}
        := \iota^* (\id - (\partial_1 - \boldsymbol{\partial}_1)h_0)^{-1}
        \quad
        \textrm{and}
        \quad
        \boldsymbol{h}_l :=
        h_l (h_{l-1} \boldsymbol{\partial}_l + \boldsymbol{\partial}_{l+1}
        h_l)^{-1}.
    \end{equation}
    Clearly the used inverse operators exist as formal power series thanks
    to (\ref{eq:ch}). The proof of the following lemma is completely
    analogously to the proofs of
    \cite[Prop.~25 and 26]{bordemann.herbig.waldmann:2000}. The
    $G$-invariance is obvious.
    \begin{lemma}
        \label{lemma:defhomotopy}
        The operators $\boldsymbol{\iota^*}$ and $\boldsymbol{h}$ are
        $G$-invariant and fulfill the relations
        \begin{equation}
            \label{eq:againhomotopy}
            \boldsymbol{h}_{l-1} \boldsymbol{\partial}_l
            +
            \boldsymbol{\partial}_{l+1} \boldsymbol{h}_l
            =
            \id_{\bigwedge^l (\mathfrak{t}) \otimes C^\infty (U)[[\lambda]]}
        \end{equation}
        as well as
        \begin{equation}
            \label{eq:stillcoboundary}
            \boldsymbol{\iota^*} \boldsymbol{\partial}_1 = 0
            \quad
            \textrm{and}
            \quad
            \boldsymbol{\iota^*} \prol =
            \id_{C^\infty(\mathcal{O})[[\lambda]]}.
        \end{equation}
    \end{lemma}
    Having the deformed restriction map and the chain homotopy it is quite
    easy to characterize the ideal generated by the `constraints' $J$:
    \begin{lemma}
        \label{lemma:keriotaimpartial}
        Let $\boldsymbol{\mathcal{I}}(J)$ be the (automatically two-sided)
        ideal generated by $J_1, \ldots, J_k$. Then the map
        $\boldsymbol{\iota^*}: C^\infty(U)[[\lambda]]
        \to C^\infty(\mathcal{O})[[\lambda]]$
        is surjective and
        \begin{equation}
            \label{eq:vanishingideal}
            \ker \boldsymbol{\iota^*}
            = \image \boldsymbol{\partial}_1
            = \boldsymbol{\mathcal{I}}(J).
        \end{equation}
    \end{lemma}
    Thus we can simply \emph{define} $\starO$ by (\ref{eq:DefRestHom})
    which gives a well-defined star-product on the quotient. It is an
    easy computation that the first order commutator of $\starO$ gives
    indeed the desired Poisson bracket. Moreover, since the $J$'s are
    real the ideal generated by them is automatically a
    $^*$-ideal. Since $h_0$ as well as $\partial$ and
    $\boldsymbol{\partial}$ are real operators, it follows that
    $\boldsymbol{\iota^*}$ is real, too.

    It remains to show that $\boldsymbol{\iota^*}$ can be written by
    use of a series of differential operators $S_r$. This is not
    completely obvious as we used the non-local homotopy $h_0$ in
    order to define $\boldsymbol{\iota^*}$. However, one can show the
    existence of the $S_r$ in the same manner as in
    \cite[Lem.~27]{bordemann.herbig.waldmann:2000}. Note that this is
    not even necessary for Corollary~\ref{corollary:posTrace}.
\end{proof}

Note that in the above construction one does not need the `full'
machinery of the BRST reduction but only the \emph{Koszul part}. The
reason is that in this case the coadjoint orbit plays the role of the
`constraint surface' \emph{and} the reduced phase space at once.

\begin{remark}
    It seems that the above statement is not the most general one 
    can obtain: There are certainly more general orbits and also 
    non-compact groups where one can find such deformed 
    restriction maps. We leave this as an open question for future 
    projects.
\end{remark}

%
%

\section{GNS representation of the positive traces}
\label{sec:gns}

Throughout this section we shall assume that $G$ is compact and
$\iota: \mathcal{O} \hookrightarrow \mathfrak{g}^*$ is a regular
orbit. Then we shall investigate the GNS representation of the
positive trace $\boldsymbol{\tau_{\mathcal{O}}}$ as constructed in the
last section.

Let us briefly recall the basic steps of the GNS construction, see
\cite{bordemann.waldmann:1998}. Having a
$^*$-algebra $\mathcal{A}$ over $\mathbb{C}[[\lambda]]$ with a
positive linear functional
$\omega: \mathcal{A} \to \mathbb{C}[[\lambda]]$ one finds that
$\mathcal{J}_\omega = \{A \in \mathcal{A} \; | \; \omega(A^*A) = 0\}$
is a left ideal of $\mathcal{A}$, the so-called Gel'fand ideal of
$\omega$. Then
$\mathfrak{H}_\omega := \mathcal{A} \big/ \mathcal{J}_\omega$ becomes
a pre-Hilbert space over $\mathbb{C}[[\lambda]]$ via
$\SP{\psi_A, \psi_B}_\omega := \omega(A^*B)$, where
$\psi_A \in \mathfrak{H}_\omega$ denotes the equivalence class of
$A$. Finally, the left representation
$\pi_\omega(A)\psi_B = \psi_{AB}$ of $\mathcal{A}$ on
$\mathfrak{H}_\omega$ turns out to be a $^*$-representation, i.e. one
has
$\SP{\psi_B, \pi_\omega(A)\psi_C}_\omega =
\SP{\pi_\omega(A)\psi_B, \psi_C}_\omega$.

According to Theorem~\ref{theorem:iotaHomo} we have in our case a
surjective $^*$-homomorphism
\begin{equation}
    \label{eq:iotaAgain}
    \boldsymbol{\iota^*}:
    C^\infty (\mathfrak{g}^*)[[\lambda]]
    \to
    C^\infty (\mathcal{O})[[\lambda]]
\end{equation}
and a positive linear functional $\boldsymbol{\tau_{\mathcal{O}}}$
which is the pull back of a positive linear functional on
$C^\infty (\mathcal{O})[[\lambda]]$ under $\boldsymbol{\iota^*}$,
namely the trace $\tr_{\mathcal{O}}$ on $\mathcal{O}$. Thus we can use
the functoriality properties of the GNS construction, see
\cite[Prop.~5.1 and Cor.~5.2]{bordemann.neumaier.waldmann:1999} in
order to relate the GNS construction for
$\boldsymbol{\tau_{\mathcal{O}}}$ with the one for
$\tr_{\mathcal{O}}$, which is well-known, see
\cite[Sect.~5]{waldmann:2000} and
\cite[Lem.~4.3]{bordemann.roemer.waldmann:1998}. Since
$\tr_{\mathcal{O}}$ is a faithful functional the GNS representation of
$C^\infty (\mathcal{O})[[\lambda]]$ with respect to
$\tr_{\mathcal{O}}$ is simply given by left multiplication
$\mathsf{L}$ with respect to $\starO$, where
$\mathfrak{H}_{\tr_{\mathcal{O}}} = C^\infty(\mathcal{O})[[\lambda]]$.
Thus we arrive at the following theorem which can also be checked
directly:
\begin{theorem}
    \label{theorem:GNSRep}
    Let $G$ be compact,
    $\iota: \mathcal{O} \hookrightarrow \mathfrak{g}^*$ a regular
    orbit, and
    $\boldsymbol{\tau_{\mathcal{O}}}
    = \tr_{\mathcal{O}} \circ \boldsymbol{\iota^*}$
    the positive trace as in (\ref{eq:posTrace}).
    \begin{enumerate}
    \item $\supp \boldsymbol{\tau_{\mathcal{O}}} =
        \iota(\mathcal{O})$.
    \item The Gel'fand ideal
        $\mathcal{J}_{\boldsymbol{\tau_{\mathcal{O}}}}$ of
        $\boldsymbol{\tau_{\mathcal{O}}}$ coincides with
        $\ker \boldsymbol{\iota^*}$.
    \item The GNS pre-Hilbert space
        $\mathfrak{H}_{\boldsymbol{\tau_{\mathcal{O}}}}$ is unitarily
        isomorphic to $C^\infty (\mathcal{O})[[\lambda]]$ endowed with
        the inner product
        $\SP{\phi, \chi}_{\mathcal{O}}
        := \tr_{\mathcal{O}} (\cc\phi \starO \chi)$
        via
        \begin{equation}
            \label{eq:Unitary}
            U: \mathfrak{H}_{\boldsymbol{\tau_{\mathcal{O}}}}
            \ni \psi_f \mapsto \boldsymbol{\iota^*} f \in
            C^\infty (\mathcal{O})[[\lambda]]
        \end{equation}
        with inverse $U^{-1}: \phi \mapsto \psi_{\prol \phi}$.
    \item For the GNS representation
        $\pi_{\boldsymbol{\tau_{\mathcal{O}}}}$ one obtains
        \begin{equation}
            \label{eq:NiceGNSRep}
            \pi_{\mathcal{O}} (f) \phi
            :=
            U \pi_{\boldsymbol{\tau_{\mathcal{O}}}} (f) U^{-1} \phi
            =
            \boldsymbol{\iota^*} ( f \starBCH \prol \phi)
            =
            \mathsf{L}_{\boldsymbol{\iota^*}f} \phi.
        \end{equation}
    \end{enumerate}
\end{theorem}

Since the group $G$ acts on $\mathcal{O}$ and since all relevant maps
are $G$-invariant/equivariant we arrive at the following
$G$-invariance of the representation. This can be checked either
directly or follows again from
\cite[Prop.~5.1 and Cor.~5.2]{bordemann.neumaier.waldmann:1999}.
\begin{lemma}
    \label{lemma:GinvRep}
    The GNS representation $\pi_{\mathcal{O}}$ is $G$-equivariant in
    the sense that
    \begin{equation}
        \label{eq:GequivRep}
        \pi_{\mathcal{O}} (g^* f) g^* \phi =
        g^* (\pi_{\mathcal{O}} (f) \phi)
    \end{equation}
    for all $\phi \in C^\infty (\mathcal{O})[[\lambda]]$,
    $f \in C^\infty (\mathfrak{g}^*)[[\lambda]]$ and
    $g \in G$. Moreover, the $G$-representation on
    $C^\infty (\mathcal{O})[[\lambda]]$ is unitary.
\end{lemma}

Let us finally mention a few properties of the commutant of
$\pi_{\mathcal{O}}$ and the `baby-version' of the Tomita-Takesaki
theory arising from this representation. The following statements
follow almost directly form the considerations in
\cite[Sect.~7]{waldmann:2000}. We consider the anti-linear map
\begin{equation}
    \label{eq:ModularConj}
    J: \phi \mapsto \cc\phi,
\end{equation}
where $\phi \in C^\infty (\mathcal{O})[[\lambda]]$, which is clearly
anti-unitary with respect to the inner product
$\SP{\cdot,\cdot}_{\mathcal{O}}$ and involutive. This map plays the
role of the \emph{modular conjugation}. The \emph{modular operator}
$\Delta$ is just the identity map since in our case the linear
functional is a trace, i.e. a \emph{KMS functional} for inverse
temperature $\beta = 0$. Then we can characterize the commutant of the
representation $\pi_{\mathcal{O}}$ as follows:
\begin{proposition}
    \label{proposition:TomitaTakesaki}
    For $f \in C^\infty (\mathfrak{g}^*)[[\lambda]]$ we denote by
    $\mathsf{R}_{\boldsymbol{\iota^*} f}$ the right multiplication
    with $\boldsymbol{\iota^*} f$ with respect to the star-product
    $\starO$. Then the map
    \begin{equation}
        \label{eq:ModConj}
        \pi_{\mathcal{O}} (f)
        = \mathsf{L}_{\boldsymbol{\iota^*} f}
        \mapsto
        J \mathsf{L}_{\boldsymbol{\iota^*} f} J
        = \mathsf{R}_{\boldsymbol{\iota^*} \cc f}
    \end{equation}
    is an anti-linear bijection onto the commutant
    $\pi_{\mathcal{O}}'$ of $\pi_{\mathcal{O}}$.
\end{proposition}
Note that in this particularly simple case the modular one-parameter
group $U_t$ is just the identity
$U_t = \id_{C^\infty (\mathcal{O})[[\lambda]]}$, since we have a
trace. More generally, one could also consider KMS functionals of the
form
$f \mapsto
\boldsymbol{\tau_{\mathcal{O}}} (\mathrm{Exp}(-\beta H) \starBCH f)$
where $H \in C^\infty (\mathfrak{g}^*)[[\lambda]]$ and $\mathrm{Exp}$
denotes the star exponential with respect to $\starBCH$ and
$\beta \in \mathbb{R}$ is the `inverse temperature'.

From the above proposition we immediately have the following result on
the relation between the $\mathfrak{g}$-representations on
$C^\infty(\mathcal{O})[[\lambda]]$ arising from the GNS construction.
\begin{lemma}
    \label{lemma:gRep}
    For $x, y \in \mathfrak{g}$ we have
    \begin{align}
        \label{eq:gLeft}
        \pi_{\mathcal{O}} (\hat x)
        \pi_{\mathcal{O}} (\hat y)
        -
        \pi_{\mathcal{O}} (\hat y)
        \pi_{\mathcal{O}} (\hat x)
        &=
        \im\lambda \pi_{\mathcal{O}} (\widehat{[x,y]})
        \\
        \label{eq:gRight}
        \mathsf{R}_{\boldsymbol{\iota^*} \hat x}
        \mathsf{R}_{\boldsymbol{\iota^*} \hat y}
        -
        \mathsf{R}_{\boldsymbol{\iota^*} \hat y}
        \mathsf{R}_{\boldsymbol{\iota^*} \hat x}
        &=
        - \im \lambda \mathsf{R}_{\boldsymbol{\iota^*} \widehat{[x,y]}}
    \end{align}
    and
    \begin{equation}
        \label{eq:gAd}
        \pi_{\mathcal{O}} (\hat x)
        - \mathsf{R}_{\boldsymbol{\iota^*} \hat x}
        =
        \im\lambda \Lie_{x_{\mathcal{O}}},
    \end{equation}
    where $\Lie_{x_{\mathcal{O}}}$ denotes the Lie derivative in
    direction of the fundamental vector field of $x$.
\end{lemma}

%
%

\section{Traces for deformations via group actions}
\label{sec:action}

Let us now describe a quite general mechanism for constructing
deformations and traces via group actions. We first consider the
algebraic part of the construction. Let $G$ be a group and denote the
right translations by $\Right_g : h \mapsto hg$, where 
$g, h \in G$. The left translations are denoted by $\Left_g$,
respectively. Moreover, let $\mathcal{A}_G \subseteq \mathrm{Fun}(G)$
be a sub-algebra of the complex-valued functions on $G$, closed under
complex conjugation. We require
$\Right_g^* \mathcal{A}_G \subseteq \mathcal{A}_G$ for all $g \in
G$. Then an associative formal deformation 
$(\mathcal{A}_G[[\lambda]],\star_G)$ of $\mathcal{A}_G$ is called
\emph{(right) universal deformation} if it is right-invariant, i.e.
\begin{equation}
    \label{eq:unidef}
    \Right_g^* (f_1 \star_G f_2) 
    = \Right_g^* f_1 \star_G \Right_g^* f_2
\end{equation}
for all $g \in G$ and $f_1, f_2 \in \mathcal{A}_G[[\lambda]]$. Thus
the right translations act as automorphisms of $\star_G$. In the
sequel we shall always assume that $1 \in \mathcal{A}_G$ and 
$1 \star_G f = f = f \star_G 1$. 
\begin{remark}
   If $G$ is a Lie group and $\mathcal{A}_G$ are all smooth functions
   on $G$ then the existence of a right-invariant deformation gives 
   quite strong conditions on $G$. However, in typical examples one 
   may only deform a smaller class of functions. For instance the data 
   of a $G$-invariant star product on a homogeneous symplectic space 
   $G\stackrel{\pi}{\to}H\backslash G$ determines a right deformation 
   of $\mathcal{A}_G:=\pi^{*}C^\infty(H\backslash G)$. In the extreme case 
   where $H=\{e\}$, the pair $(\mathcal{A}_G,\star_G)$ becomes a star 
   product algebra $(C^\infty(G)[[\lambda]],\star_\lambda)$. 
   The Poisson structure on $G$ associated to the first order term of 
   $\star_\lambda$ is then right-invariant. Its characteristic 
   distribution (generated by Hamiltonian vector fields)---being 
   integrable and right-invariant--- determines a Lie subalgebra 
   $\mathcal{S}$ of $\mathfrak{g} = \mbox{Lie}(G)$ endowed with a 
   non-degenerate Chevalley 2-cocycle $\Omega$ with respect to the trivial 
   representation of $\mathcal{S}$ on $\mathbb{R}$. This type of Lie 
   algebras $(\mathcal{S},\Omega)$ (or rather their associated Lie groups) has 
   been studied by Lichnerowicz et al.. When unimodular such a Lie
   algebra is solvable \cite{lichnerowicz.medina:1988}.
\end{remark}

Now consider a set $X$ with a left action $\tau: G \times X \to X$ of
$G$. For abbreviation we shall sometimes write $g.x$ instead of
$\tau(g,x)$. We shall use the universal deformation $\star_G$ in order
to induce a deformation of a certain sub-algebra of
$\mathrm{Fun}(X)$. First we define 
$\alpha^x: \mathrm{Fun}(X) \to \mathrm{Fun}(G)$ by
\begin{equation}
    \label{eq:alphaxdef}
    (\alpha^x f)(g) = (\tau^*_g f)(x)
\end{equation}
for $x \in X$ and $g \in G$. Having specified $\mathcal{A}_G$ we
define the space
\begin{equation}
    \label{eq:AXdef}
    \mathcal{A}_X 
    = \{ f \in \mathrm{Fun}(X) 
    \; | \; \alpha^x f \in \mathcal{A}_G 
    \quad \textrm{for all} \quad x \in X \},
\end{equation}
which is clearly a sub-algebra of $\mathrm{Fun}(X)$ stable under
complex conjugation. Let us remark that $\mathcal{A}_X$ contains at
least those functions on $X$ which are constant along the orbits of
$\tau$. Indeed, let $f \in \mathrm{Fun}(X)$ satisfy $f(g.x) = f(x)$
for all $x \in X$ and $g \in G$. Then 
$(\alpha^x f)(g) = f(g.x) = f(x)$ is constant (not depending on $g$).

The deformation $\star_G$ induces canonically an associative
deformation $\star_X$ of $\mathcal{A}_X$, thereby justifying the name
`universal deformation'. Indeed, define
\begin{equation}
    \label{eq:starX}
    (f_1 \star_X f_2)(x) = (\alpha^x f_1 \star_G \alpha^x f_2) (e),
\end{equation}
where $e \in G$ denotes the unit element. Then we have the following
proposition:
\begin{proposition}
    \label{proposition:starX}
    Let $(\mathcal{A}_G[[\lambda]], \star_G)$ be a universal
    deformation and $(\mathcal{A}_X[[\lambda]], \star_X)$ as above.
    \begin{enumerate}
    \item Then $(\mathcal{A}_X[[\lambda]], \star_X)$ is an associative
        formal deformation of $\mathcal{A}_X$ which is Hermitian if
        $\star_G$ is Hermitian. Moreover, 
        $\alpha^x: 
        (\mathcal{A}_X[[\lambda]], \star_X) \to
        (\mathcal{A}_G[[\lambda]], \star_G)$ 
        is a homomorphism of associative algebras.
    \item If $f_1$ is constant on some orbit $G.x_0$ then 
        \begin{equation}
            \label{eq:constorbit}
            (f_1 \star_X f_2)(g.x_0) = f_1(g.x_0) f_2(g.x_0)
            = (f_2 \star_X f_1)(g.x_0)
        \end{equation}
        for all functions $f_2 \in \mathcal{A}_X[[\lambda]]$.
        In particular, the $\star_X$-product with a function, which
        is constant along all orbits, is the undeformed product. Thus
        $\star_X$ is `tangential' to the orbits in a very strong
        sense.
    \end{enumerate}
\end{proposition}
\begin{proof}
    Let us first recall a few basic properties of $\alpha^x$, $\tau$,
    $\Right$, and $\Left$. The following relations are straightforward
    computations:
    \begin{equation}
        \label{eq:Ralpha}
        \Right_g^* \alpha^x = \alpha^{g.x}
        \quad \textrm{and} \quad
        \Left_g^* \alpha^x = \alpha^x \tau_g^*.
    \end{equation}
    Using the right invariance of $\star_G$ and the above rules we
    find the following relation
    \begin{equation}
        \label{eq:homomorphism}
        \alpha^x (f_1 \star_X f_2) = \alpha^x f_1 \star_G \alpha^x f_2
    \end{equation}
    for $f_1, f_2 \in \mathcal{A}_X[[\lambda]]$. This implies on one
    hand that $\mathcal{A}_X[[\lambda]]$ is indeed closed under the
    multiplication law $\star_X$. On the other hand it follows that
    $\alpha^x$ is a homomorphism. With \eqref{eq:homomorphism} the
    associativity of $\star_X$ is a straightforward
    computation. Finally, if $\star_G$ is Hermitian then $\star_X$ is
    Hermitian, too, since all involved maps are real, i.e. commute
    with complex conjugation. For the second part one computes
    \begin{equation}
        \label{eq:constorb}
        (f_1 \star_X f_2) (g.x_0) 
        = (\alpha^{x_0} f_1 \star_G \alpha^{x_0} f_2) (g).
    \end{equation}
    Now $\alpha^{x_0} f_1$ is constant whence the $\star_G$-product is
    the pointwise product. Thus the claim easily follows. If this
    holds even for all orbits and not just for $G.x_0$ then the
    $\star_X$-product with $f_1$ is the pointwise product globally.
\end{proof}
\begin{remark}
    \label{remark:starXnontrivial}
    From \eqref{eq:constorbit} we conclude that, heuristically
    speaking, the deformation $\star_X$ becomes more non-trivial the
    larger the orbits of $\tau$ are.
\end{remark}
\begin{remark}
    \label{remark:Leftuniversal}
    Given a right universal deformation $(\mathcal{A}_G, \star_R)$,
    one gets a left universal deformation $(\mathcal{A}_G, \star_L)$
    via the formula 
    \begin{equation}
        \label{eq:leftuniv}
        a \star_L b = \iota^*(\iota^* a \star_R \iota^* b)
    \end{equation}
    provided $\mathcal{A}_G$ is a bi-invariant subspace. Here 
    $\iota: G \to G$ denotes the inversion map $g \to g^{-1}$.
    Starting with a left invariant deformation 
    $(\mathcal{A}_G,\star_G)$ of $G$ and an action 
    $\tau: G \times X \to X$, the associated deformation of
    $\mathcal{A}_X$ is then defined by the formula 
    \begin{equation}
        \label{eq:leftinduceddef}
        (f_1 \star f_2) (x)
        = (\iota^* \alpha^x f_1 \star_G \iota^* \alpha^x f_2)(e).
    \end{equation}
\end{remark}
In some interesting cases, in particular in the Abelian case, the
universal deformation $\star_G$ is also \emph{left invariant},
i.e. the left translations $\Left^*_g$ acts as automorphisms of
$\star_G$, too. In this situation the induced deformation $\star_X$ is
invariant under $\tau_g^*$:
\begin{lemma}
    \label{lemma:tauinv}
    Let $\mathcal{A}_G$ be in addition left invariant and let
    $\star_G$ be a bi-invariant universal deformation. Then
    $\mathcal{A}_X$ is invariant under $\tau_g^*$ for all $g \in G$
    and
    \begin{equation}
        \label{eq:invstarX}
        \tau_g^* (f_1 \star_X f_2) = \tau_g^* f_1 \star_X \tau^*_g f_2.
    \end{equation}
\end{lemma}
\begin{proof}
    This is a straightforward computation using only the definitions
    and \eqref{eq:Ralpha}.
\end{proof}

Our main interest in the universal deformations comes from the
following simple observation:
\begin{theorem}
    \label{theorem:main}
    Let $(\mathcal{A}_G[[\lambda]], \star_G)$ be a right universal
    deformation and let 
    $\tr_G: \mathcal{A}_G[[\lambda]] \to \mathbb{C}[[\lambda]]$ be a
    trace with respect to $\star_G$. Let 
    $\Phi: \mathrm{Fun}(X)[[\lambda]] \to \mathbb{C}[[\lambda]]$ be an
    arbitrary $\mathbb{C}[[\lambda]]$-linear functional. Then
    $\tr_\Phi: \mathcal{A}_X[[\lambda]] \to \mathbb{C}[[\lambda]]$
    defined by
    \begin{equation}
        \label{eq:tracePhi}
        \tr_\Phi (f) = \Phi( x \mapsto \tr_G (\alpha^x f))
    \end{equation}
    is a trace with respect to $\star_X$.
\end{theorem}
\begin{proof}
    This follows directly from the homomorphism property of $\alpha^x$
    and the trace property of $\tr_G$.
\end{proof}

In particular the trace $\tr_G$ combined with the evaluation
functionals at some point $x \in X$
\begin{equation}
    \label{eq:trx}
    \tr_x: f \mapsto \tr_G (\alpha^x f)
\end{equation}
yields a trace for $\star_X$. Thus the only difficult task is to find
traces for $\star_G$.

As a last remark we shall discuss the positivity of the traces
$\tr_\Phi$. We assume that $\tr_G$ is a positive trace whence
$\tr_G(\cc{f} \star_G f) \ge 0$ in the sense of formal power series
for all $f \in \mathcal{A}_G[[\lambda]]$.
\begin{lemma}
    \label{lemma:pos}
    Assume $\tr_G$ is a positive trace and $\Phi$ takes non-negative
    values on non-negative valued functions on $X$. Then $\tr_\Phi$ is
    positive. In particular $\tr_x$ is always positive.
\end{lemma}
\begin{remark}
    \label{remark:convergent}
    The above construction has the big advantage that it can be
    transfered to the framework of topological deformations instead of
    formal deformations. This has indeed been done by Rieffel
    \cite{rieffel:1993} in a $C^*$-algebraic framework for actions of
    $\mathbb{R}^d$. For a class of non-abelian groups this has been 
    done in \cite{bieliavsky.massar:2000:pre}.
\end{remark}

Let us finally mention two examples. The first one is the well-known
example of the Weyl-Moyal product for $\mathbb{R}^{2n}$ and the second
is obtained as the asymptotic version of
\cite{bieliavsky.massar:2000:pre} for rank one Iwasawa subgroups of
$\mathrm{SU}(1,n)$.
\begin{example}
    \label{example:weyl}
    Let $\starW$ be the Weyl-Moyal star product on $\mathbb{R}^{2n}$,
    explicitly given by
    \begin{equation}
        \label{eq:weylagain}
        f \starWeyl g = 
        \mu \circ \eu^{\frac{\lambda}{2\im}\sum_k
          (\partial_{q^k} \otimes \partial_{p_k} 
          - \partial_{p_k} \otimes \partial_{q^k})} 
        f \otimes g,
    \end{equation}
    where $\mu(f \otimes g) = fg$ is the pointwise product and 
    $q^1 \ldots, p_n$ are the canonical Darboux coordinates on
    $\mathbb{R}^{2n}$. Clearly $\starWeyl$ is invariant under
    translations whence it is a bi-invariant universal deformation of
    $C^\infty(\mathbb{R}^{2n})[[\lambda]]$. Moreover, it is well-known
    that $\starWeyl$ is strongly closed, whence the integration with 
    respect to
    the Liouville measure provides a trace, which is positive. Thus
    one can apply the above general results to this situation.
\end{example}
\begin{example}
    \label{example:iwasawa}
    This example is the asymptotic version of
    \cite{bieliavsky.massar:2000:pre}. 
    The groups we consider are Iwasawa subgroups $G=AN$ of
    $\mathrm{SU}(1,n)$, where $\mathrm{SU}(1,n)=ANK$ is an Iwasawa
    decomposition. One has the obvious $G$-equivariant diffeomorphism 
    $G \to \mathrm{SU}(1,n)/K$ (here $K = \mathrm{U}(n)$). The group
    $G$ therefore inherits a left-invariant symplectic (K{\"a}hler)
    structure coming from the one on the rank one Hermitian symmetric
    space $\mathrm{SU}(1,n)/\mathrm{U}(n)$. The symplectic group may
    then be described as follows. As a manifold, one has 
    \begin{equation}
        \label{eq:Gdiffeo}
        G = \mathbb{R} \times \mathbb{R}^{2n} \times \mathbb{R}.
    \end{equation}
    In these coordinates the group multiplication law reads
    \begin{equation}
        \label{eq:groupmult}
        \Left_{(a,x,z)}(a',x',z')
        =\left(
            a+a', 
            \eu^{-a'} x + x', 
            \eu^{-2a'} z + z' + \frac{1}{2}\Omega(x,x') \eu^{-a'}
        \right),
    \end{equation}
    where $\Omega$ is a constant symplectic structure on the vector
    space $\mathbb{R}^{2n}$. The 2-form 
    \begin{equation}
        \label{eq:sympstructG}
        \omega = \Omega + da \wedge dz
    \end{equation}
    then defines a left-invariant symplectic structure on $G$.
    The universal deformation $\starBM$ we are looking for is a star
    product for this symplectic structure. Since on
    $\mathbb{R}^{2n+2}$ all symplectic star products are equivalent,
    it will be sufficient to describe $\starBM$ be means of an
    equivalence transformation 
    $T = \id + \sum_{r=1}^\infty \lambda^r T_r$ relating $\starBM$ and
    $\starWeyl$. In \cite{bieliavsky.massar:2000:pre} an explicit
    integral formula for $T$ has been given, which is defined on the
    Schwartz space $\mathcal{S}(\mathbb{R}^{2n+2})$. It allows for an
    asymptotic expansion in $\hbar$ and gives indeed the desired
    equivalence transformation $T$.
    Then $\starBM$ defined by
    \begin{equation}
        \label{eq:starBM}
        f \starBM g = T^{-1} (Tf \starWeyl Tg)
    \end{equation}
    is a left-invariant universal deformation of $G$ and again we can
    use this to apply the above results on universal deformations. 
    Moreover, since $\starWeyl$ is strongly closed, the functional
    \begin{equation}
        \label{eq:traceforstarBM}
        \tr^G (f) := \int_G T(f) \; \omega^{n+1}
    \end{equation}
    defines a trace functional for $\starBM$ on
    $C^\infty_0(G)[[\lambda]]$. This is again positive since that $T$ is real
    i.e. $\cc{Tf}=T\cc{f}$. 
    
    In what follows we give a precise 
    description of the star product $\starBM$ in the two dimensional 
    case i.e. on the group $ax+b$. The higher dimensional case is 
    similar but more intricate. The non-formal deformed product in 
    the $ax+b$ case is obtained by transforming Weyl's product on 
    $(\mathbb{R}^2,da\wedge d\ell)$ under the equivalence
    \begin{equation}
        T = F^{-1}\circ\phi_\hbar^*\circ F
    \end{equation}
    where
    \begin{equation}
        Fu(a,\alpha) = \int \eu^{-\im\alpha \ell} \, u(a,\ell) 
        \; d\ell 
        \quad
        \textrm{with}
        \quad 
        u \in \mathcal{S}(\mathbb{R}^2)
    \end{equation}
    is the partial Fourier transform in the second variable and where 
    $\phi_\hbar:\mathbb{R}^2\to\mathbb{R}^2$ is the one-parameter family 
    of diffeomorphisms given by 
    \begin{equation}
        \phi_\hbar(a,\alpha)
        = (a,\frac{1}{\hbar}\sinh(\alpha\hbar))
        \quad(\hbar\in\mathbb{R}).
    \end{equation}
    One has 
    \begin{eqnarray}
    Tu(a, \ell) = c \int \eu^{\im\alpha \ell} 
    \eu^{\frac{-\im}{\hbar}\sinh(\alpha\hbar)q} \, u(a,q) 
    \; dq\,d\alpha
    = c \int \eu^{\im\alpha(\ell-q)} 
    \eu^{-\im\psi_\hbar(\alpha)q} \,u(a,q)
    \; dq\,d\alpha
    \end{eqnarray}
    with
    \begin{equation}
    \psi_\hbar(\alpha)
    = \sum_{k\geq 1}\frac{\hbar^{2k}\alpha^{2k+1}}{(2k+1)!}.
    \end{equation}
    Setting $p=\hbar\alpha$, one gets
    \begin{equation}
    Tu(a,\ell)=\frac{c}{\hbar}
    \int \eu^{\frac{\im}{\hbar}p(\ell-q)}
    \eu^{\frac{-\im}{\hbar}\psi_1(p)q} \, u(a,q) \; dq\,dp
    \end{equation}
    which precisely coincides with
    \begin{equation}
        \label{NFEQ}
        \id\otimes Op_{\hbar,1}(\eu^{\frac{-\im}{\hbar}\psi_1(p)q})) 
        \, u(a,\ell)
    \end{equation}
    where $Op_{\hbar,1}f(p,q)$ denotes the anti-normally ordered 
    quantization of the function $f(q,p)$. Recall that the  
    $\kappa$-ordered 
    pseudodifferential quantization rule on $(\mathbb{R}^2,dq\wedge dp)$ is 
    defined (at the level of test functions) by
    $Op_{\hbar,\kappa}:\mathcal{D}(\mathbb{R}^2)
    \longrightarrow\End(L^2(\mathbb{R}))$ with
    \begin{equation}
    Op_{\hbar,\kappa}(f)\varphi(q)=\frac{c}{\hbar} \int 
    \eu^{\frac{\im}{\hbar}p(q-\xi)} \, 
    f\left(\kappa\xi+(1-\kappa)q,p\right)\varphi(\xi)
    \; d\xi\,dp\quad(\kappa\in[0,1]).
    \end{equation}
    The explicit asymptotic expansion formula for $Op_{\hbar,\kappa}(f)$ 
    is well known, see e.g. 
    \cite[Sect.~1.2, p.~231 and Eq.~(58), p.~258]{stein:1993}.
    It yields an expression for the equivalence $T$
    at the formal level which we write, with natural delicacy, as
    \begin{equation}
        T = \id\otimes\exp\left(\frac{\im}{\lambda}
        \psi_{1}(\frac{\lambda}{\im}\partial_\ell).\ell\right),
    \end{equation}
    where the operator 
    $T_{(\ell)}:=\exp\left(\frac{\im}{\lambda}
    \psi_1(\frac{\lambda}{\im}\partial_\ell).\ell\right)$ 
    is to be understood as 
    anti-normally ordered ($\kappa=1$). Observe the reality of the 
    equivalence, which
may be directly checked using the fact that the function $\psi_1$ is 
odd. Moreover, for every right-invariant vector field 
$X$ on $G=ax+b$, one checks \cite{bieliavsky.massar:2000:pre} that 
$T\circ X\circ T^{-1}$ is an inner 
derivation of the Moyal-Weyl product $\starWeyl$. In other words, the 
star product $\starBM$ is left-invariant on G. 
\end{example}

%
%


\begin{thebibliography}{10}

\bibitem {arnal.benamar:2000}
{\sc Arnal, D., Ben~Amar, N.: }\newblock {\em {K}ontsevich's Wheels and
  Invariant Polynomial Functions on the Dual of Lie Algebras}.
\newblock Lett. Math. Phys.  {\bf 52} (2000), 291--300.

\bibitem {arnal.benamar.masmoudi:1999}
{\sc Arnal, D., Ben~Amar, N., Masmoudi, M.: }\newblock {\em Cohomology of Good
  Graphs and {K}ontsevich Linear Star Products}.
\newblock Lett. Math. Phys.  {\bf 48} (1999), 291--306.

\bibitem {basart.et.al:1984}
{\sc Basart, H., Flato, M., Lichnerowicz, A., Sternheimer, D.: }\newblock {\em
  Deformation Theory applied to Quantization and Statistical Mechanics}.
\newblock Lett. Math. Phys.  {\bf 8} (1984), 483--394.

\bibitem {bayen.et.al:1978}
{\sc Bayen, F., Flato, M., Fr{{\o}}nsdal, C., Lichnerowicz, A., Sternheimer,
  D.: }\newblock {\em Deformation Theory and Quantization}.
\newblock Ann. Phys.  {\bf 111} (1978), 61--151.

\bibitem {bertelson.cahen.gutt:1997}
{\sc Bertelson, M., Cahen, M., Gutt, S.: }\newblock {\em Equivalence of Star
  Products}.
\newblock Class. Quantum Grav.  {\bf 14} (1997), A93--A107.

\bibitem {bieliavsky.massar:2000:pre}
{\sc Bieliavsky, P., Massar, M.: }\newblock {\em Strict Deformation
  Quantizations for Actions of a Class of Symplectic Lie Groups}.
\newblock Preprint  {\bf math.QA/0011144} (November 2000).

\bibitem {bordemann.herbig.waldmann:2000}
{\sc Bordemann, M., Herbig, H.-C., Waldmann, S.: }\newblock {\em BRST
  Cohomology and Phase Space Reduction in Deformation Quantization}.
\newblock Commun. Math. Phys.  {\bf 210} (2000), 107--144.

\bibitem {bordemann.neumaier.waldmann:1998}
{\sc Bordemann, M., Neumaier, N., Waldmann, S.: }\newblock {\em Homogeneous
  Fedosov Star Products on Cotangent Bundles I: Weyl and Standard Ordering with
  Differential Operator Representation}.
\newblock Commun. Math. Phys.  {\bf 198} (1998), 363--396.

\bibitem {bordemann.neumaier.waldmann:1999}
{\sc Bordemann, M., Neumaier, N., Waldmann, S.: }\newblock {\em Homogeneous
  Fedosov star products on cotangent bundles II: GNS representations, the WKB
  expansion, traces, and applications}.
\newblock J. Geom. Phys.  {\bf 29} (1999), 199--234.

\bibitem {bordemann.roemer.waldmann:1998}
{\sc Bordemann, M., R{\"{o}}mer, H., Waldmann, S.: }\newblock {\em A Remark on
  Formal KMS States in Deformation Quantization}.
\newblock Lett. Math. Phys.  {\bf 45} (1998), 49--61.

\bibitem {bordemann.waldmann:1998}
{\sc Bordemann, M., Waldmann, S.: }\newblock {\em Formal GNS Construction and
  States in Deformation Quantization}.
\newblock Commun. Math. Phys.  {\bf 195} (1998), 549--583.

\bibitem {bursztyn.waldmann:2000b}
{\sc Bursztyn, H., Waldmann, S.: }\newblock {\em Deformation Quantization of
  Hermitian Vector Bundles}.
\newblock Lett. Math. Phys.  {\bf 53} (2000), 349--365.


\bibitem {bursztyn.waldmann:2000a}
{\sc Bursztyn, H., Waldmann, S.: }\newblock {\em On Positive Deformations of
  {$^*$}-Algebras}.
\newblock In: {\sc Dito, G., Sternheimer, D. (eds.): }\newblock {\em
  Conf{\`e}rence Mosh{\`e} Flato 1999. Quantization, Deformations, and
  Symmetries}, {\em Mathematical Physics Studies} no. {\bf 22},   69--80.
  Kluwer Academic Publishers, Dordrecht, Boston, London, 2000.

\bibitem {bursztyn.waldmann:2001a}
{\sc Bursztyn, H., Waldmann, S.: }\newblock {\em Algebraic Rieffel Induction,
  Formal Morita Equivalence and Applications to Deformation Quantization}.
\newblock J. Geom. Phys.  {\bf 37} (2001), 307--364.

\bibitem {cahen.gutt.rawnsley:1996}
{\sc Cahen, M., Gutt, S., Rawnsley, J.: }\newblock {\em On Tangential Star
  Products for the Coadjoint Poisson Structure}.
\newblock Commun. Math. Phys.  {\bf 180} (1996), 99--108.

\bibitem {connes:1994}
{\sc Connes, A.: }\newblock {\em Noncommutative Geometry}.
\newblock Academic Press, San Diego, New York, London, 1994.


\bibitem {connes.flato.sternheimer:1992}
{\sc Connes, A., Flato, M., Sternheimer, D.: }\newblock {\em Closed Star
  Products and Cyclic Cohomology}.
\newblock Lett. Math. Phys.  {\bf 24} (1992), 1--12.

\bibitem {dito:1999}
{\sc Dito, G.: }\newblock {\em {K}ontsevich Star Product on the Dual of a Lie
  Algebra}.
\newblock Lett. Math. Phys.  {\bf 48} (1999), 307--322.

\bibitem{dito.sternheimer:2002}
{\sc Dito, G., Sternheimer, D.: }\newblock
{\em Deformation Quantization: Genesis, Developments and
  Metamorphoses}.
\newblock To appear in the Proceedings of the meeting between
mathematicians and theoretical physicists, Strasbourg, 2001. IRMA
Lectures in Math. Theoret. Phys., vol. 1, Walter De Gruyter, Berlin
2002, pp. 9--54.

\bibitem {fedosov:1996}
{\sc Fedosov, B.~V.: }\newblock {\em Deformation Quantization and Index
  Theory}.
\newblock Akademie Verlag, Berlin, 1996.

\bibitem {felder.shoikhet:2000}
{\sc Felder, G., Shoikhet, B.: }\newblock {\em Deformation Quantization with
  Traces}.
\newblock Lett. Math. Phys.  {\bf 53} (2000), 75--86.

\bibitem {giaquinto.zhang:1998}
{\sc Giaquinto, A., Zhang, J.~J.: }\newblock {\em Bialgebra actions, twists,
  and universal deformation formulas}.
\newblock J. Pure Appl. Algebra  {\bf 128}.2 (1998), 133--152.

\bibitem {gutt:1983}
{\sc Gutt, S.: }\newblock {\em An Explicit $*$-Product on the Cotangent Bundle
  of a Lie Group}.
\newblock Lett. Math. Phys.  {\bf 7} (1983), 249--258.

\bibitem {gutt:2000}
{\sc Gutt, S.: }\newblock {\em Variations on deformation quantization}.
\newblock In: {\sc Dito, G., Sternheimer, D. (eds.): }\newblock {\em
  Conf{\`e}rence Mosh{\`e} Flato 1999. Quantization, Deformations, and
  Symmetries}, {\em Mathematical Physics Studies} no. {\bf 21},   217--254.
  Kluwer Academic Publishers, Dordrecht, Boston, London, 2000.

\bibitem {gutt.rawnsley:1999}
{\sc Gutt, S., Rawnsley, J.: }\newblock {\em Equivalence of star products on a
  symplectic manifold; an introduction to Deligne's {\v{C}}ech cohomology
  classes}.
\newblock J. Geom. Phys.  {\bf 29} (1999), 347--392.

\bibitem {gutt.rawnsley:2001:pre}
{\sc Gutt, S., Rawnsley, J.: }\newblock {\em Traces for star products on
  symplectic manifolds}.
\newblock Preprint  {\bf math.QA/0105089} (May 2001).

\bibitem {karabegov:1998b}
{\sc Karabegov, A.~V.: }\newblock {\em On the Canonical Normalization of a
  Trace Density of Deformation Quantization}.
\newblock Lett. Math. Phys.  {\bf 45} (1998), 217--228.

\bibitem {kontsevich:1997:pre}
{\sc Kontsevich, M.: }\newblock {\em Deformation Quantization of Poisson
  Manifolds, I}.
\newblock Preprint  {\bf q-alg/9709040} (September 1997).

\bibitem {kontsevich:1999}
{\sc Kontsevich, M.: }\newblock {\em Operads and Motives in Deformation
  Quantization}.
\newblock Lett. Math. Phys.  {\bf 48} (1999), 35--72.

\bibitem{lichnerowicz.medina:1988} {\sc Lichnerowicz, A., Medina, A.:} 
\newblock {\em Groupes a structures symplectiques ou kaehleriennes invariantes.} 
\newblock C. R. Acad. Sci., Paris, Ser. I 306, No.3  (1988), 133--138.

\bibitem {nest.tsygan:1995a}
{\sc Nest, R., Tsygan, B.: }\newblock {\em Algebraic Index Theorem}.
\newblock Commun. Math. Phys.  {\bf 172} (1995), 223--262.

\bibitem {omori.maeda.yoshioka:1991}
{\sc Omori, H., Maeda, Y., Yoshioka, A.: }\newblock {\em Weyl Manifolds and
  Deformation Quantization}.
\newblock Adv. Math.  {\bf 85} (1991), 224--255.

\bibitem {rieffel:1993}
{\sc Rieffel, M.~A.: }\newblock {\em Deformation quantization for actions of
  $\mathbb R^d$}.
\newblock Mem. Amer. Math. Soc.  {\bf 106}.506 (1993).

\bibitem {rosenberg:1996:pre}
{\sc Rosenberg, J.: }\newblock {\em Rigidity of K-theory under deformation
  quantization}.
\newblock Preprint  {\bf q-alg/9607021} (July 1996).

\bibitem {sternheimer:1998}
{\sc Sternheimer, D.: }\newblock {\em Deformation Quantization: Twenty Years
  After}. \newblock
In: {\sc Rembieli\`{n}ski, J. (ed.):} \newblock
{\em Particles, Fields, and Gravitation}.\newblock
AIP Press, New York 1998.

\bibitem{stein:1993}
{\sc Stein, E.M.: }\newblock {\em Harmonic Analysis Real-Variable Methods, Orthogonality, \& Oscillatory
Integrals}.
\newblock Princeton Mathematical Series, Princeton University Press (1993).

\bibitem{tamarkin.tsygan:2001}
{\sc Tamarkin, D., Tsygan, B.: }\newblock
{\em Cyclic Formality and Index Theorems}.
\newblock Lett. Math. Phys. {\bf 56} (2001), 85--97.

\bibitem {waldmann:2000}
{\sc Waldmann, S.: }\newblock {\em Locality in GNS Representations of
  Deformation Quantization}.
\newblock Commun. Math. Phys.  {\bf 210} (2000), 467--495.

\bibitem {weinstein:1994}
{\sc Weinstein, A.: }\newblock {\em Deformation Quantization}.
\newblock S\'eminaire Bourbaki 46\`eme ann\'ee  {\bf 789} (1994).

\bibitem {weinstein:1997}
{\sc Weinstein, A.: }\newblock {\em The modular automorphism group of a Poisson
  manifold}.
\newblock J. Geom. Phys.  {\bf 23} (1997), 379--394.

\bibitem {weinstein.xu:1998}
{\sc Weinstein, A., Xu, P.: }\newblock {\em Hochschild cohomology and
  characterisic classes for star-products}.
\newblock In: {\sc Khovanskij, A., Varchenko, A., Vassiliev, V. (eds.):
  }\newblock {\em Geometry of differential equations. Dedicated to V. I. Arnold
  on the occasion of his 60th birthday},   177--194. American Mathematical
  Society, Providence, 1998.

\end{thebibliography}

\end{document}